\documentclass[11pt]{article}\usepackage{graphicx}\usepackage[T1]{fontenc}
\usepackage[utf8]{inputenc}
\usepackage[dvips]{epsfig}
\usepackage{amssymb}
\newtheorem{theorem}{Theorem} \newtheorem{lemma}{Lemma}\newtheorem{proposition}{Proposition}\newtheorem{claim}{Claim}
\newtheorem{corollary}{Corollary} \newcommand{\La}{\Lambda}
\newcommand{\R}{{\mathbb R}}  \newcommand{\Z}{{\mathbb Z}} \newcommand{\N}{{\mathbb N}}
\newcommand{\T}{{\mathbb T}} \newcommand{\C}{{\mathbb C}}  \newcommand{\D}{{\mathbb D}}

\begin{document}\title{On Szeg\"{o}--Kolmogorov Prediction Theorem}

\author{Alexander Olevskii  and  Alexander Ulanovskii}

\date{}\maketitle

\section{Introduction}

        Let $w$ be a weight on the circle $\T=[0,1] $ (a nonnegative function
        with finite integral).
        Given a discrete set $\La\subset\R$, we denote by  $E(\La)$ the family
        of exponentials
$$
                      E(\La): = \{e^{2\pi i \lambda t} , \lambda\in \La\}.
$$

        The classical Szeg\"{o}--Kolmogorov
        theorem (see \cite{GS58, K41})  implies that the family $E(\N),   \N :=\{1,2,3,...\}, $
       spans the whole weighted space $L^2(\T,w)$ if and only if
\begin{equation}\label{int}
\int_{0}^{1}\log w(t)\,dt=-\infty.
\end{equation}

        Kolmogorov's probabilistic interpretation of this result connects it with the
      possibility to  `predict precisely the future from the past' for the  stationary
      stochastic processes  with discrete time. Namely, the prediction is possible
      if and only if the density $w$ of the spectral measure of the process satisfies
      condition (\ref{int}).  In particular, this is true when $w$ has an
     `exponentially deep zero' at some point.
      There is a number of papers of different authors related to this famous result.

        In the present note we are interested in the following problem:
       Can the `prediction' be possible if some part of the `past' is not known?
       And how large this part can be?
        In other words, we wonder how many exponentials with positive frequencies can be removed
      while still keeping the completeness property in the same weighted spaces.
          Problems of this type, in a more general form,  were posed by A. N. Kolmogorov,       see \cite{K45}, referring to Kolmogorov's letter.

      We consider the problem under the restriction that a single `deep zero' is the
       only reason for divergence of the logarithmic integral (\ref{int}). Our main results show
       that in this situation the frequencies
       $\{n^3: n\in\N\}$ can be removed with no affect on the completeness,
        while $\{n^2: n\in\N\}$ can't.

      We will consider the class of weights $w$ which have a `deep   right-hand zero' at the origin:

        \medskip\noindent
(A) $w$ is positive, increasing and bounded on $(0,1)$ and satisfies (\ref{int}).


\medskip
 Note that the results below remain true if $w$ is only increasing on some subinterval $(0,\delta)\subset (0,1)$.
 Note also that one could equally well consider  a `deep   right-hand zero' at any other point $t_0\in\T$.

   Here are the precise statements of our results:

            \begin{theorem}\label{t1}  The system $E(\N\setminus\Gamma), \Gamma\subset \N,$ spans $L^2(\T,w) $ for every weight $w$  satisfying (A),
                         whenever
\begin{equation}\label{gamma}
\sum_{\gamma\in \Gamma}\frac{1}{\sqrt\gamma}< \infty.
\end{equation}
\end{theorem}

\begin{theorem}\label{t2} The system $E(\N\setminus\{k^2:k\in\N\})$ does not span $L^2(\T,w)$, with some $w$ satisfying (A).
\end{theorem}

     More generally, an analogue of Theorem \ref{t2}  remains true for every regularly
             distributed  $\Gamma$  with  infinite sum in (\ref{gamma})        and  with an appropriate weight $w $ satisfying (A),
           see sec. \ref{s5} below.




           Observe that 
           assumption (A) in Theorem \ref{t1} is essential.
          Using some constructions of periodic functions with `small'
          support and spectrum, one can see that for the general weights $w$ satisfying (\ref{int})
          and  sets $\Gamma$ satisfying condition (\ref{gamma})  (and even a much stronger sparseness condition
          of this type),  the system $E(\N \setminus\Gamma)$ may be incomplete in $L^2(\T,w)$ with some $w$ satisfying (\ref{int}), see sec. \ref{s6}.


\section{Extension  to Non-Harmonic Series}\label{s2}

Given a function $H(t)$  defined on $\R$, we write $H\in L^2(\T,1/w)$ if $$\int_0^1|H(t)|^2\frac{1}{w(t)}\,dt<\infty.$$

 Assume  that $w$ satisfies  (A). Then $1/w(t)>C,t\in\T,$ for some positive constant $C$, and so $L^2(\T,1/w)\subset L^2(\T)$.

   Recall that the system   $E(\N\setminus\Gamma)$ is not complete in $L^2(\T,w)$ if and only if there is a function  $\Phi\in L^2(\T,w)$  orthogonal to every element of $E(\N\setminus\Gamma)$:
  $$
  \int_0^1 \Phi(t)e^{-2\pi i n t}w(t)\,dt=0,\quad n\in\N\setminus\Gamma.
  $$
   Clearly, the function $H(t):=\Phi(t)w(t)$ belongs to $L^2(\T,1/w)\subset L^2(\T)$,
  and so its Fourier coefficients vanish on $\N\setminus\Gamma$. Hence, $E(\N\setminus\Gamma)$ is not complete in $L^2(\T,w)$ if and only if there is a non-trivial Fourier series \begin{equation}\label{h00}H(t)=\sum_{n\in\Z_-}a_ne^{2\pi i nt}+\sum_{n\in\Gamma}b_ne^{2\pi i n t},\quad \{a_n\},\{b_n\}\in  l^2,\end{equation}
which belongs to  $L^2(\T,1/w)$.
We see that Theorem \ref{t1} is equivalent to the statement that if $\Gamma\subset\N$ satisfies (\ref{gamma}) then no non-trivial function $H$ in (\ref{h00}) may belong to $L^2(\T,1/w)$, for every $w$ satisfying (A).

We  extend Theorem \ref{t1} to non-harmonic series
\begin{equation}\label{h01}
H(t)=\sum_{\lambda\in\La}a_\lambda e^{2\pi i \lambda t}+\sum_{\gamma\in\Gamma}b_\gamma e^{2\pi i \gamma t},\quad \{a_\lambda\},\ \{b_\gamma\}\in l^2,
\end{equation}
where $\La\subset(-\infty,0]$ and $\Gamma\subset(0,\infty)$ are arbitrary uniformly discrete sets.

Recall that a discrete set $\La\subset\R$ is called uniformly discrete (u.d.) if
$$
\inf_{\lambda,\lambda'\in\La,\lambda\ne\lambda'}|\lambda-\lambda'|>0.
$$
If both $\La$ and $\Gamma$ are u.d., it is well-known that every function $H$ in (\ref{h01}) belongs to $L^2(I)$ on every finite interval $I\subset\R$.

\begin{theorem}\label{t4} Assume  $\La\subset[-\infty,0)$ and $\Gamma\subset(0,\infty)$ are u.d. sets. If $\Gamma$ satisfies  (\ref{gamma}) then no non-trivial function $H$ in (\ref{h01}) belongs  to $L^2(\T,1/w)$, for every $w$ satisfying (A).
\end{theorem}

Clearly,  it also follows that for every $t_0\in\R$ we have  $H(t+t_0)\not\in L^2(\T,1/w),$ for every $w$ satisfying (A).
\section{Proofs of Theorems \ref{t1} and \ref{t4}}

Recall that, as explained in sec. \ref{s2},  Theorem \ref{t1} follows from Theorem \ref{t4}.

Throughout the rest of this section  we assume that $\Gamma\subset(0,\infty)$ is a u.d. set satisfying  (\ref{gamma}).

Our proof of Theorem \ref{t4} relies on a construction of a function $g$ satisfying the properties:
\begin{lemma}\label{l0}There is a function $g\in \mathcal{C}^1(\R)$ such that
\begin{equation}\label{g0}g(t)=0,\ \ t\geq 0,\quad  g, g'\in L^1(\R),\end{equation} and whose Cauchy transform vanishes on $\Gamma$:
\begin{equation}\label{g1}\int_{-\infty}^0\frac{g(u)}{u-\gamma}\,du=0,\quad \gamma\in\Gamma.\end{equation}
\end{lemma}

Lemma \ref{l0} will be proved at the end of  next subsection.

\subsection{Auxiliary Lemmas}
Denote by $\C_+ := \{z \in\C : \Im z > 0\}$ the (open) upper half plane, and by $\C_-$  the lower one. As usual, by $\mathcal{H}^p(\C_+)$  and  $\mathcal{H}^p(\C_-)$ we mean  the Hardy spaces in the upper and lower half-plane.

  A fundamental
fact from the Hardy spaces theory states: Assume $f\in \mathcal{H}^p(\C_-),1\leq p\leq\infty,$ then the (non-tangential) boundary values of $f$ exist a.e. and satisfy
$$\int_{-\infty}^\infty\frac{\log|f(t)|}{1+t^2}\,dt>-\infty.$$

Below we will use an easy consequence of this fact:

\begin{corollary}\label{c1} Assume a nontrivial function $f(z)$ is analytic in $\C_-$ and satisfies $$|f(z)|\leq C +
C|z|,\quad z \in \C_-,$$ for some $C > 0$. Then $f$ has (non-tangential) boundary values a.e. on $\R$ satisfying
$$\int_0^{1}\log|f(t)|\,dt>-\infty.$$\end{corollary}

Consider the domain
$D := \C\setminus (-\infty,0]$. We fix the (principal) value of argument, $-\pi<\arg z<\pi, z\in D,$  and set
\begin{equation}\label{varphi}\varphi(z) :=
\frac{z^2}{(1+\sqrt z)^8}\prod_{\gamma\in\Gamma}\left(1-\sqrt{\frac{z}{\gamma}}\right)\left(1+\sqrt{\frac{z}{\gamma}}\right)^{-1},\quad z\in D.
\end{equation}

The product in the right hand-side of (\ref{varphi}) converges
 due to (\ref{gamma}). It is easy to check that it  is a Blaschke product  in $D$ (i.e. it is analytic in $D$, its modulus is less than one in $D$ and  is equal to one on the boundary of $D$). So, $\varphi$ is analytic in $D$. Clearly, $\varphi$ vanishes on $\Gamma$ and satisfies$$\varphi(z) =\overline{\varphi(\bar z)},\quad z\in D.$$

The boundary of $D$  consists of two copies of $(-\infty, 0]$. Define the boundary values of $\varphi$ by
$\varphi_{\pm}(t), t\in  (-\infty, 0]$, where $\varphi_+(t)$ is the limit of $\varphi(z)$ as $z$ tends to $t$ from `above'. It is easy to
check that $\varphi$ is bounded in $D$, continuous up to the boundary and
\begin{equation}\label{varphi1}
\varphi_+(t)=\overline{\varphi_-(t)} =
\frac{-t^2}{(1+i\sqrt{|t|})^8}\prod_{\gamma\in\Gamma}
\left( 1-i\sqrt{\frac{|t|}{\gamma}}\right)
\left(1+i\sqrt{\frac{|t|}{\gamma}}\right)^{-1},  t\leq 0.
\end{equation}

Denote by ${\bf 1}_{(-\infty,0]}(t)$ the characteristic function of $(-\infty,0]$.

Our first lemma follows easily from (\ref{varphi1}) by elementary calculations:

\begin{lemma}\label{l1} The function $\varphi_+(t)\cdot{\bf 1}_{(-\infty,0]}(t),t\in\R,$ is continuously differentiable and satisfies
$$\frac{d}{dt}\left(\varphi_+(t)\cdot{\bf 1}_{(-\infty,0]}(t)\right)\in L^2(\R).$$\end{lemma}

The  the next lemma is well-known (and simple) and is given without proof:

\begin{lemma}\label{l2} Assume a continuously differentiable function $f(t)$ satisfies $f,f'\in L^2(\R)$. Then its
Fourier transform belongs to $L^1(\R)$.\end{lemma}

We now introduce the function $g$ which appears in Lemma \ref{l0}: $$g(t):= (\varphi_+(t)-\varphi_-(t))\cdot {\bf 1}_{(-\infty,0)}(t)=2i\Im\varphi_+(t)\cdot{\bf 1}_{(-\infty,0)}(t).$$
Then $g(t)=0, t\geq0$ and we have\begin{equation}\label{v0}
\int_{\partial D}\frac{\varphi(t)}{t-z}\,dt=\int_{-\infty}^0\frac{g(t)}{t-z}\,dt,\quad z\in\D,
\end{equation}where the contour integration along $\partial D$ is in the positive direction with respect to  $D$.

\begin{lemma}\label{l3} We have\begin{equation}\label{h0}|g(t)|\leq\frac{2t^2}{(1+t)^4},\quad t\leq0,\quad
\mbox{and}\quad g'\in L^2(\R);\end{equation}
\begin{equation}\label{varphi2}
\varphi(z)=\frac{1}{2\pi i}\int_{-\infty}^0\frac{g(t)}{t-z}\,dt,\quad z\in D.\end{equation}
\end{lemma}

Indeed, the absolute value of the (Blaschke) product in (\ref{varphi}) is $1$ on $\partial D$,  so that $$|\varphi(t)|= \frac{t^2}{|1+\sqrt t|^8}=\frac{t^2}{(1+|t|)^4}, \quad t\leq0.$$
 This proves the first condition in (\ref{h0}).  Also, clearly, $$|\varphi(z)|=O(|z|^{-2}),\quad z\in D,|z|\to\infty.$$Then  representation (\ref{varphi2})  follows from (\ref{v0}) and
the elementary residue theory. The second condition in (\ref{h0}) follows from Lemma \ref{l1}. 


\begin{lemma}\label{l5} Given two locally $L^2$-functions $F$ and $\Phi$ on $\R$ which vanish on $(-\infty,0)$. If $F\in L^2(\T,1/w)$ with some $w$ satisfying (A), then
 $F \ast  \Phi\in L^2(\T,1/w)$.\end{lemma}

It follows from (A) that  $1/w$ is decreasing on  $(0,1)$. Hence,
$$
\int_0^1 |(F\ast G)(t)|^2\frac{1}{w(t)}\,dt=\int_0^1 \frac{1}{w(t)}\left|\int_0^t G(t-s)F(s)\,ds\right|^2dt\leq
$$$$ \|G\|^2_{L^2(0,1)}\int_0^1\int_0^t |F(s)|^2\frac{1}{w(s)}\,dsdt<\infty.
$$

\begin{lemma}\label{l6} Assume $\Phi\in L^2(\T,1/w)$ with some $w$ satisfying (\ref{int}). Then
$$\int_0^1\log| \Phi(t)|\,dt=-\infty.$$\end{lemma}

Indeed, since the exponential function $e^x$ is  convex,
 Jensen's inequality yields
 $$
 \exp\left(\int_0^1 \left(2\log|\Phi(t)|-\log w(t)\right)\,dt\right)\leq \int_0^1|\Phi(t)|^2\frac{1}{w(t)}\,dt= \|\Phi\|^2_{L^2(\T,1/w)}<\infty.
 $$Now, the lemma follows from (\ref{int}).

Finally, observe that  Lemma \ref{l0} stated above is an easy consequence of the definitions of $\varphi$ and $g$ and Lemma \ref{l3}.

\subsection{Proof of Theorem \ref{t4}}\label{s3}

Throughout the proof we denote by C absolute positive constants.
We may assume that $\Gamma$ contains an infinite number of elements.

Assume a function $H$ in (\ref{h01}) satisfies $H\in L^2(\T,1/w),$ with some $w$ satisfying (A).  To prove Theorem \ref{t4}, we have to show that $H=0$ a.e.

The proof  consists of several steps.

1. Let $g$ be the function in Lemma \ref{l0}, and denote by $G$ its inverse Fourier transform:
$$G(s): =
\int_{-\infty}^{0}e^{2\pi i st}g(t)\,dt,\quad s\in\R.$$
Since $g,g'\in L^1(\R)$, by Lemma \ref{l2}  we have $G\in L^1(\R)$.

2. For every real $a$ set$$
e_a(t) := e^{2\pi i at}\cdot {\bf 1}_{[0,\infty)}(t).$$

Set
$$H_\La(t) :=\sum_{\lambda\in\La}a_\lambda e_\lambda(t),\quad H_\Gamma(t):=\sum_{\gamma\in\Gamma}b_\gamma e_\gamma(t).$$
Then $$H_+(t) := H(t)\cdot{\bf 1}_{(0,\infty)}(t)=H_\La(t)+H_\Gamma(t).$$ Clearly, to prove Theorem \ref{t4} it suffices to prove that $H_+=0$ a.e.

3. Denote $G_+(t):=G(t)\cdot{\bf 1}_{(0,\infty)}(t)$, and consider the convolution
$$H_+\ast G_+ =H_\La\ast G_++H_\Gamma\ast G_+.$$
The main step of the proof is to show that both convolutions $H_\La\ast G_+$ and $H_\Gamma\ast G_+$  admit analytic
continuation from $(0,\infty)$ to the lower half-plane $\C_-$. We will show that the continuation of the first
one belongs to $H^\infty(\C_-)$, while the absolute value of the second one is less than $C +C|z|, z\in \C_-$.
Since $H_+\in L^2(\T,1/w)$, by Lemma~\ref{l5}, $H_+\ast G_+\in L^2(\T,1/w)$. Hence, using Corollary~\ref{c1} and Lemma \ref{l6}, we
may conclude that $H_+\ast G_+ = 0$, and so $H_+ = 0$.

4. Given a real number $q$, consider the convolution
$$(e_q \ast G_+)(t) =\int_0^t e^{2\pi qi (t-s)}G(s)\,ds=\int_0^t e^{2\pi i q(t-s)}\int_{-\infty}^0e^{2\pi i us}g(u)\,du\,ds$$
\begin{equation}\label{eq8}
=\frac{1}{2\pi i}\int_{-\infty}^0\frac{e^{2\pi i qt}-e^{2\pi i ut}}{u-q}g(u)\,du,\quad t>0.
\end{equation}

Assume $q\in \Gamma$. Then by (\ref{g1}),
$$(e_q \ast  G_+)(t) =\frac{1}{2\pi i}\int_{-\infty}^0\frac{e^{2\pi i ut}}{q-u}g(u)\,du,\quad t>0, \ q\in\Gamma.$$ This gives
\begin{equation}\label{eq9}
(H_\Gamma\ast  G_+)(t) =\frac{1}{2\pi i}\int_{-\infty}^0e^{2\pi i ut}\left(\sum_{\gamma\in\Gamma}\frac{b_\gamma}{\gamma-u}g(u)\right)\,du,\quad t>0.\end{equation}

Clearly, since $\Gamma$ is u.d. and lies on the positive half-axis, for every $u \leq 0$ we have
$$
\sum_{\gamma\in\Gamma}\left|\frac{b_\gamma}{\gamma-u}\right|\leq\left(\sum_{\gamma\in\Gamma}|b_\gamma|^2\right)^{1/2}
\left(\sum_{\gamma\in\Gamma}\frac{1}{\gamma^2}\right)^{1/2}< C < 1.$$
Since $g\in L^1(\R)$,  the integral in the right side of (\ref{eq9}) admits analytic continuation
to $\C_-$ and is bounded there.

5. Assume $q\leq 0$. For every fixed $u\ne q, u < 0,$ the entire function
$$\frac{e^{2\pi i q z}-e^{2\pi i u z}}{2\pi iz (u-q)}$$
is the inverse Fourier transform of the function ${\bf 1}_{(u,q)}(t)/(q-u)$. It is easy to check that the estimate is true:
$$\left|\frac{e^{2\pi i qz }-e^{2\pi i u z}}{2\pi i (u-q)}\right|\leq |z|\max\{e^{-2\pi qy},e^{-2\pi u y}\}\leq |z|,\quad z = x + iy, y< 0.$$

Since $\La$ is u.d., for every fixed $u \leq 0$ there exists $n\in\N$, such that at most $n$ elements $\lambda\in\La$
satisfy $|\lambda-u|\leq 1.$ Denoting by $K$ the maximal sum of $|a_\lambda|$ of $n$ consequent elements of $\La$, we get
$$
\sum_{\lambda\in\La}\left|a_\lambda\frac{e^{2\pi i \lambda z}-e^{2\pi i uz}}{2\pi i(u-\lambda)}\right|\leq $$$$ K|z|+\left(\sum_{\lambda\in\La}|a_\lambda|^2\right)^{1/2}
\left(\sum_{\lambda\in\La, |\lambda-u|>1}\frac{1}{|u-\lambda|^2}\right)^{1/2}\leq K|z|+C,\quad z\in \C_-.$$
By (\ref{eq8}), we see that the convolution $H_\La\ast G_+$ admits analytic continuation to $\C_-$ and is $\leq C +C|z|$
there. This finishes the proof of Theorem \ref{t4}.

\section{Proof of Theorem \ref{t2}}\label{s4}

As explained in sec. 2, to prove Theorem \ref{t2} it suffices to find a nontrivial Fourier series $H$  whose Fourier coefficients vanish on $\N\setminus \{k^2,k\in\N\}$,
and $H\in L^2(\T,1/w)$ for some $w$ satisfying (A). We show that 
$H$ can be chosen with a symmetric spectrum concentrated on $\pm (2k+1)^2$:

\begin{theorem}\label{t20}
There exists a non-trivial Fourier series $H$, $$H(t)=\sum_{k\in\N}c_k\sin(2\pi (2k+1)^2t),$$satisfying
\begin{equation}\label{we}
|H(t)|\leq C_1e^{-C_2/|t|},\quad |t|<1,
\end{equation}where $C_1$ and $C_2$ are some positive constants.
\end{theorem}

Clearly, the function $H$ in Theorem \ref{t20} belongs to $L^2(\T,1/w),$ with   the weight $w(t)=\exp(-2C_2/t), t\in\T,$ satisfying condition (A).

Condition (\ref{we}) means that $H$ has a `deep zero' at the origin.
 It should be mentioned that connection between the `deepness of zero' of
    a function and the degree of lacunarity of its Fourier series was studied
    by  Mandelbrot \cite{m}.  See \cite{l}, Appendix 2,   for further results and     references.

Set
\begin{equation}\label{ggg}
h(z):=e^{i\pi z}\frac{\sin(\pi z)}{\cos\left(\frac{\pi}{2}\sqrt z\right)\cos\left(\frac{\pi}{2}\sqrt{- z}\right)}.
\end{equation}

 Denote by $C$ different positive constants.
 
We will need

\begin{lemma}\label{ll}
The estimates holds
\begin{equation}\label{0}|h(x+iy)|\leq \frac{Ce^{-C\sqrt y}}{1+x^2}, \quad x+iy\in\C, y\geq0.\end{equation}
\begin{equation}\label{01}|h(x+iy)|\leq \frac{Ce^{2\pi |y|-C\sqrt{ |y|}}}{1+x^2}, \quad x+iy\in\C, y<0.\end{equation}
\end{lemma}

The proof is straightforward.

It follows from (\ref{0})
that  $\{h(k),k\in\Z\}\in l^2$. Set $$H(t):=2i\sum_{k=1}^{\infty}h((2k+1)^2)\sin(2\pi(2k+1)^2t).$$

\begin{lemma}\label{ls} Denote by $\hat h(t)$ the Fourier transform of $h$. Then

(i) $\hat h\in L^2(\R)$ and $\hat h(t)=0$ a.e. for $t$ outside of the interval $(0,1);$

(ii)  $h(n)=0$ for $n\in \Z\setminus\{\pm (2k+1)^2:k\in\N\};$

(iii) $h(t)=H(t), 0<t<1;$

(iv) $H$ satisfies  (\ref{we}).
\end{lemma}

Theorem \ref{t20} is an immediate consequence of statement (iv).

Condition (i) follows from  (\ref{0}) and  (\ref{01}) and the Paley--Wiener theorem.

Condition (ii) easily follows from (\ref{ggg}).

Recall that the $k$-th Fourier coefficient of $\hat h$ is equal to $h(k)$. Clearly,  $h(-k)=-\overline{h(k)}$. Hence, condition (iii)  follows from (ii) and the definition of $H$.

Let us prove  (\ref{we}) for $0<t<1$. For every $y>0$ we  write
$$
H(t)=\hat h(t)=\int_{-\infty}^\infty e^{-2\pi i tx}h(x)\,dx=\int_{-\infty}^\infty e^{-2\pi i t(x+iy)}h(x+iy)\,dx.
$$
Hence, by (\ref{0}),
$$
|H(t)|\leq Ce^{2\pi ty-C\sqrt y}\int_{-\infty}^{\infty}\frac{dx}{1+x^2}=Ce^{2\pi ty-C\sqrt y}.
$$
Choose $y=(C/4\pi t)^{2}$, where $C$ is the constant in the exponent above, to get
$$
|H(t)|\leq Ce^{-C^2/(8\pi t)},\quad 0<t<1.
$$

 Now, use a similar argument and (\ref{01}) to check that$$|H(t)|=|\hat h(t)|\leq Ce^{-C^2/(8\pi(1- t))},\quad 0<t<1.$$ 
Due to periodicity of  $H$, this proves  (\ref{we}) for $-1<t<0$.

\section{Necessity of (\ref{gamma}) for Regularly Distributed $\Gamma$}\label{s5}

As mentioned in sec. 1, the sufficiency condition for completeness (\ref{gamma}) in Theorem \ref{t1} is also necessary for a subclass of `regularly distributed' sets $\Gamma\subset\N $.

Assume $\Gamma\subset\N$ satisfies \begin{equation}\label{gamma1}
\sum_{\gamma\in \Gamma}\frac{1}{\sqrt\gamma}= \infty.
\end{equation}
 Consider the condition

\medskip
(B)  $\Gamma$ contains a subset $\Gamma_1$  which satisfies  (\ref{gamma1}) and the condition
 \begin{equation}\label{gamma2} \#(\Gamma_1\cap(0,x))\leq C\frac{\sqrt x}{\log x}(1+o(1)),\quad x\to\infty,\end{equation}
 where $C>0$ is sufficiently small fixed number.

\begin{theorem}\label{tt} Assume $\Gamma\subset\N$  satisfies (B).
 Then   $E(\Z\setminus\Gamma )$ is not complete in $L^2(\T,w)$ for some   $w$ satisfying (A).
\end{theorem}

We omit the proof of Theorem \ref{tt}, which is somewhat similar to the proof of Theorem \ref{t2} but is much more technical.

\medskip\noindent
{\bf Remark}.
1. Condition (B) means that one may make $\Gamma$ sparser to satisfy condition (\ref{gamma2}) while preserving condition (\ref{gamma1}).

\medskip\noindent 2. Here we present an example which  tests the sharpness of condition (\ref{gamma}).
Set $$\Gamma_\rho:=\{[k^\rho], k\in\N\},$$where $[\cdot]$ means the integer part. One may check that $\Gamma_\rho$  satisfies (B) for $\rho\leq2$, and that it satisfies (\ref{gamma}) for $\rho>2$. Hence, Theorems \ref{t1} and \ref{tt}  show that   $E(\N\setminus\Gamma_\rho)$ is complete in $L^2(\T,w),$ for every $w$ satisfying (A) if and only if $\rho>2$.


\medskip\noindent 3.  If a weight $w$ has a   `very deep zero'  then most part of $\N$ can be removed
  without loss of completeness. In particular, one can remove certain sets $\Gamma$ of lower density one:

\begin{proposition}\label{p1}   There is a set $\Gamma\subset\N$ of lower density one,
such that $E(\N\setminus\Gamma)$ is complete in $L^2(\T,w)$ with some  positive increasing weight $w$ satisfying (\ref{int}).
 \end{proposition}


For the weights $w_0$ vanishing on some subinterval $(0,\delta)\subset\T$, this can       be deduced from
 the classical  Beurling--Malliavin theorem \cite{BM67}.
         By an appropriate    perturbation of $w_0$, one can  construct a weight $w$ as in Proposition \ref{p1}.

\section{Functions with Small Support and Spectrum}\label{s6}

        In the results above the divergence of logarithmic integral in (\ref{int}) is
        due to the singularity localized near a single point.         However, the divergence may also occur  when the singularity
        is `spread out' along a set of positive measure.
         In this  case the removal of even a `very sparse' set  of frequencies may lead to incompleteness of remaining exponentials
        with respect to some weight satisfying (\ref{int}).

         The example below if  based on the following

   \begin{claim}
     There is a function $H \in L^2(\T)$ satisfying the properties

           (i)   $H = 0$ on a set $E\subset\T$ of positive measure;

           (ii)  The spectrum of $H$ is sparse (in the sense clarified below).\end{claim}

           There are several approaches to construction of such functions $\varphi$ which employ  different
           meanings of  the term `sparseness'.

             An approach we use here was inspired by the classical Men'shov Correction Theorem,  saying that every measurable function $f : \T\to \C$ can be changed
           on a (dense) set  of arbitrarily small measure $ \varepsilon>0$ so that the Fourier series of the `corrected' function
           $f_\epsilon$ converges uniformly.

  A  modification of the proof of this theorem shows that the spectrum of $f_\varepsilon$ can be chosen to be quite sparse (see  \cite{O85}).

             Arutunyan \cite{A84}  obtained a nice version of Menshov's theorem
            in which the `sparseness' is understood  as follows: {\sl The spectrum of $f_\varepsilon$ can be localized inside of any   given in advance symmetric set $\La\subset\Z$, which contains arbitrarily long         intervals of integers}.

         Consider the property of a set $\Gamma\subset\N:$

   \medskip
            (C) $\Gamma$ contains arbitrarily long         intervals of integers.
            \medskip

Let us prove

\begin{proposition}If $\Gamma$ satisfies (C) then the system $E(\Z\setminus (\Gamma\cup(-\Gamma))$ does not span  $L^2(\T,w)$, for some $w$ satisfying (\ref{int}).\end{proposition}

Indeed, assume $\Gamma$ satisfies (C). Choose any non-trivial continuous function $f$ which vanishes on an interval in $\T$. Applying  Arutynyan's theorem with a sufficiently small $\varepsilon$, one gets a non-trivial function
            $H:= f_\epsilon\in L^2(\T)$ with spec$\,H\subset\Gamma\cup(-\Gamma)$, which vanishes on some  set $E\subset\T$ of positive measure.

         Let $w$ be any weight satisfying (\ref{int}) and such that  $w(t)\geq 1$ on $\T\setminus E $.  Then $H \in L^2(\T,1/w)$. This implies that $E(\Z\setminus (\Gamma\cup(-\Gamma))$
         does not span
         $L^2(\T,w).$ 

           Clearly,  sets $\Gamma$ satisfying (C) can be `very sparse'. In particular, one may prove the following statement: {\sl
          For every positive no matter how slowly decreasing function  $\psi(t)\to0,t\to\infty$, there is a set $\Gamma\subset\N$ such that
          $$
          \sum_{\gamma\in\Gamma}\psi(\gamma)<\infty,
          $$and the system $E(\N\setminus\Gamma)$ does not span $L^2(\T,w)$, for some $w$ satisfying (\ref{int}).}

          \noindent A.O.: School of Mathematics, Tel Aviv University\\ Ramat Aviv,  69978 Israel. E-mail:
 olevskii@post.tau.ac.il

\medskip

\noindent A.U.: Stavanger University,  4036 Stavanger, Norway\\ E-mail: Alexander.Ulanovskii@uis.no

\end{document}